\documentclass[12pt]{amsart}
\usepackage{amsthm}
\usepackage{amssymb}
\usepackage{cite}
\usepackage{geometry}
\usepackage{enumerate}
\usepackage{setspace}
\usepackage{tikz}
\usepackage{tikz-cd}
\usetikzlibrary{matrix, arrows}
\usepackage[unicode=true]
 {hyperref}
\hypersetup{
colorlinks=true,
urlcolor=black,
citecolor=blue,
linkcolor=blue,
}

\allowdisplaybreaks

\newtheorem{thm}{Theorem}
\newtheorem{prop}[thm]{Proposition}

\theoremstyle{definition}

\theoremstyle{remark}
\newtheorem{remark}[thm]{Remark}

\numberwithin{equation}{section}

\newcommand{\Aut}{\mathrm{Aut}}
\newcommand{\Hol}{\mathrm{Hol}}

\newcommand{\Inn}{\mathrm{Inn}}
\newcommand{\PSL}{\mathrm{PSL}}
\newcommand{\PGL}{\mathrm{PGL}}

\begin{document}

\large 

\title[Almost simple groups whose holomorph contains a solvable regular subgroup]{Finite almost simple groups whose holomorph\\ contains a solvable regular subgroup}
\author{Cindy (Sin Yi) Tsang}
\address{Department of Mathematics, Ochanomizu University, 2-1-1 Otsuka, Bunkyo-ku, Tokyo, Japan}
\email{tsang.sin.yi@ocha.ac.jp}
\urladdr{http://sites.google.com/site/cindysinyitsang/} 

\date{\today}

\maketitle


\begin{abstract}In our previous paper, we gave a complete list of the finite non-abelian simple groups whose holomorph contains a solvable regular subgroup. In this paper, we  refine our previous work by considering all finite almost simple groups. In particular, our result yields a complete characterization of the finite almost simple groups which occur as  the type of a Hopf--Galois structure on a solvable extension, or equivalently, the additive group of a skew brace having a solvable multiplicative group.\end{abstract}


\tableofcontents


\section{Introduction}

Let $G$ and $N$ be two finite groups of the same order. We say that the pair $(G,N)$ is \textit{realizable} if the holomorph of $N$ contains a regular subgroup that is isomorphic to $G$. Recall that the \textit{holomorph of $N$} is defined to be 
\[ \Hol(N) = \rho(N)\rtimes\Aut(N) = \lambda(N) \rtimes \Aut(N)\]
as a subgroup of the symmetric group of $N$, where $\lambda$ and $\rho$ are the left and right regular representations of $N$. By \cite{GP,By96} and \cite{Skew braces}, respectively, it is known that $(G,N)$ is realizable if and only if there exists a Hopf-Galois structure of type $N$ on a Galois $G$-extension of fields, or equivalently, a skew brace with additive group $N$ and multiplicative group isomorphic to $G$. This is exactly our motivation of studying realizability (see \cite{book} and \cite{SB} for more details).

\vspace{1.5mm}

For the pair $(G,N)$ to be realizable, it is conjectured that if $G$ is insolvable then $N$ is also insolvable (see \cite{Tsang solvable,Byott conj} for some progress on this). On the other hand, the converse is known to be false -- there exist realizable pairs $(G,N)$ for which $G$ is solvable but $N$ is insolvable. It seems to be a natural problem to try and determine the insolvable groups $N$ that can occur in such a pair. In a previous paper \cite[Theorems 1.2 \& 1.3]{Tsang last}, we proved that:

\begin{thm}\label{old thm}Let $N$ be a finite almost simple group. If $(G,N)$ is realizable for some solvable group $G$, then the socle of $N$ is isomorphic to
\begin{enumerate}[$(1)$]
\item $\PSL_3(3), \, \PSL_3(4),\, \PSL_3(8), \, \mathrm{PSU}_3(8),\, \mathrm{PSU}_4(2),\, \mathrm{M}_{11}$; or
\item $\PSL_2(q)$ with $q\neq 2,3$ a prime power.
\end{enumerate}
The converse is also true when $N$ is a finite non-abelian simple group.
\end{thm}

The purpose of this paper is to study the converse of Theorem \ref{old thm} for an arbitrary finite almost simple group $N$ that is not necessarily simple. As it turns out, the proposition below (see \cite[Corollary 2.2 \& Proposition 2.7]{Tsang last}) is enough to do the job except when $N = \mathrm{PSU}_3(8)$. This case has already been dealt with in \cite{Tsang last} and so may be ignored. %
\begin{prop}\label{key prop}Let $N$ be any finite group.
\begin{enumerate}[$(a)$]
\item If $(G,N)$ is realizable for a solvable group $G$, then there exists a subgroup
\[ \Inn(N) \leq P \leq \Aut(N)\]
such that $P=AB$ is factorized by some solvable subgroups $A$ and $B$ that satisfy the equality $A\Inn(N) =B\Inn(N)$.
\vspace{1.5mm}
\item If $N$ has trivial center and there exists a subgroup
\[ \Inn(N) \leq P \leq \Aut(N)\]
such that $P=AB$ is factorized by some solvable subgroups $A$ and $B$ that satisfy $A\cap B=1,\, A\Inn(N) =B\Inn(N)$, and $A$ splits over $A\cap \Inn(N)$, then $(G,N)$ is realizable for some solvable group $G$.
\end{enumerate}
Here $\Inn(N)$ denotes the inner automorphism group of $N$.
\end{prop}

The converse of Theorem \ref{old thm} does not hold in general, but we have a very nice criterion and here is our main result. 

\begin{thm}\label{new thm}Let $N$ be a finite almost simple group whose socle $\mathrm{Soc}(N)$ is isomorphic to one of the non-abelian simple groups in Theorem \ref{old thm}. There is a solvable group $G$ such that $(G,N)$ is realizable if and only if
\begin{enumerate}[$(1)$]
\item $\mathrm{Soc}(N)\simeq \PSL_2(q)$ with $q\neq 2,3$ a prime power; or
\vspace{1.5mm}
\item $\mathrm{Soc}(N)\simeq \PSL_3(3),\, \PSL_3(8),\, \mathrm{PSU}_4(2),\, \mathrm{M}_{11}$; or
\vspace{1.5mm}
\item $\mathrm{Soc}(N)\simeq \PSL_3(4)$ with $[N:\mathrm{Soc}(N)]\neq 2,4$; or
\vspace{1.5mm}
\item $\mathrm{Soc}(N)\simeq \PSL_3(4)$ with $[N:\mathrm{Soc}(N)]=2$ and $N\unlhd \Aut(\mathrm{Soc}(N))$; or
\vspace{1.5mm}
\item $\mathrm{Soc}(N)\simeq \mathrm{PSU}_3(8)$ with $[N:\mathrm{Soc}(N)]\neq 2,6$; or
\vspace{1.5mm}
\item $\mathrm{Soc}(N)\simeq \mathrm{PSU}_3(8)$ with $[N:\mathrm{Soc}(N)]=6$ and $N\unlhd \Aut(\mathrm{Soc}(N))$.
\end{enumerate}
\end{thm}

\begin{remark}In the language of Hopf-Galois structures and skew braces, our result yields a complete classification of the finite almost simple groups that can occur as the type of a Hopf-Galois structure on a solvable extension, or equivalently, the additive group of a skew brace with solvable multiplicative group.
\end{remark}
 
By definition, an almost simple group $N$ corresponds to a subgroup of the outer automorphism group of its socle $\mathrm{Soc}(N)$. Let us note that
\[ \mathrm{Out}(\PSL_3(3)) \simeq C_2, \,\ \mathrm{Out}(\PSL_3(4))\simeq D_{12},\,\ \mathrm{Out}(\PSL_3(8))\simeq C_6,\]
\[ \mathrm{Out}(\mathrm{PSU}_3(8))\simeq C_3\times S_3,\,\ \mathrm{Out}(\mathrm{PSU}_4(2))\simeq C_2,\,\ \mathrm{Out}(\mathrm{M}_{11}) \simeq C_1.\]
Moreover, for a prime power $q=p^f$ with $q\neq 2,3$, we have
\[\mathrm{Out}(\PSL_2(q))\simeq C_{\gcd(2,q-1)}\times C_f.\]
Here, as usual $C_n$ is the cyclic group of order $n$, $D_{2n}$ is the dihedral group of order $2n$, and $S_n$ is the symmetric group on $n$ letters. (See \cite{Wilson} for example.) 

\vspace{1.5mm}

Therefore, the converse of Theorem  \ref{old thm} holds with only five exceptions of $N$, as follows. Here, only the isomorphism class of $N$ matters, so it suffices to consider the subgroups of $\mathrm{Out}(\mathrm{Soc}(N))$ up to conjugacy.
\begin{itemize}
\item The subgroups $\mathrm{PSL}_3(4)\leq N \leq \Aut(\PSL_3(4))$ corresponding to
\[ \langle s\rangle,\,\ \langle rs\rangle,\,\ \langle r^3,s\rangle.\]
Note that $\langle r^3\rangle$ is also a subgroup of $D_{12}$ of order $2$, but this is normal\par\noindent and so falls under case (4) of Theorem \ref{new thm}. Here $r$ and $s$ are the usual generators of $D_{12}$ with $r^6 = s^2 = 1$ and $srs=r^{-1}$.
\vspace{1.5mm}
\item The subgroups $\mathrm{PSU}_3(8)\leq N \leq \Aut(\mathrm{PSU}_3(8))$ corresponding to
\[ \{1\}\times \langle (1\, 2)\rangle\mbox{ and }C_3\times \langle (1\, 2)\rangle.\]
Note that $\{1\}\times S_3$ is also a subgroup of $C_3\times S_3$ of order $6$, but this is normal and so falls under case (6) of Theorem \ref{new thm}.
 \end{itemize}
We now move on and use Proposition \ref{key prop} to prove Theorem \ref{new thm}. We use the help of \textsc{Magma} \cite{magma} except for the infinite family $\PSL_2(q)$ with $q\neq 2,3$.

\section{Proof of cases (2) to (6)}\label{sec}

We ran the following code in \textsc{Magma} \cite{magma}.

\begin{verbatim}
T:=(the non-abelian simple group in question);
AutT:=PermutationGroup(AutomorphismGroup(T));
InnT:=Socle(AutT);
NN:=[N`subgroup:N in Subgroups(AutT)|InnT subset N`subgroup];
//A list of the finite almost simple groups N with socle T
//up to conjugacy class in Aut(T).
for N in NN do
  AutN:=Normalizer(AutT,N);
  //We identify Aut(N) with the normalizer of N in Aut(T) 
  //and in particular Inn(N) with N.
  PP:=[P`subgroup:P in Subgroups(AutN)|N subset P`subgroup];
  //A list of the subgroups P between Inn(N) and Aut(N)
  //up to conjugacy class in Aut(N).
  PPa:=[]; PPb:=[];
  //A list to contain the P that satisfies the conditions
  //in Propositions 2(a) and 2(b), respectively.
  for P in PP do
    SolSub:=[S`subgroup:S in SolvableSubgroups(P)];
    if exists{<A,B>:A,B in SolSub|
         #P eq #A*#B/#(A meet B) and
         sub<P|A,N> eq sub<P|B,N>} then
      Append(~PPa,P);
    end if;
    if exists{<A,B>:A,B in SolSub|
         #(A meet B) eq 1 and
         #P eq #A*#B and
         sub<P|A,N> eq sub<P|B,N> and
         HasComplement(A,A meet N)} then
      Append(~PPb,P);
    end if;
  end for;
  if IsEmpty(PPa) then
    conclusion:="false";
    //Proposition 2(a): If PPa is empty then (G,N) is not 
    //realizable for any solvable group G.
    //Thus the converse of Theorem 1 is "false" for N.
  elif not IsEmpty(PPb) then
    conclusion:="true";
    //Proposition 2(b): If PPb is non-empty then (G,N) is 
    //realizable for some solvable group G.
    //Thus the converse of Theorem 1 is "true" for N.
  else
    conclusion:="inconclusive";
    //Proposition 2 is not sufficient to determine if the
    //converse of Theorem 1 holds for N.
  end if;
  if IsNormal(AutT,N) then
    normalORnot:="normal";
  else 
    normalORnot:="not normal";
  end if;
  <#N/#InnT,normalORnot,conclusion>;
end for;
\end{verbatim}

For $T = \PSL_3(3),\, \PSL_3(8),\, \mathrm{PSU}_4(2),\, \mathrm{M}_{11}$, the output of \texttt{conclusion} is all \texttt{true} and this yields case (2) of Theorem \ref{new thm}. 

\vspace{1.5mm}

For $T=\PSL_3(4)$ and $[N:\mathrm{Soc}(N)]\neq 2,4$, the output of \texttt{conclusion} is all \texttt{true} and this yields case (3) of Theorem \ref{new thm}.

\vspace{1.5mm}

For $T=\PSL_3(4)$ and $[N:\mathrm{Soc}(N)]= 2,4$, the output is as follows:

\begin{verbatim}
<2, "normal", "true">
<2, "not normal", "false">
<2, "not normal", "false">
<4, "not normal", "false">
\end{verbatim}

\noindent We then deduce case (4) of Theorem \ref{new thm}.

\vspace{1.5mm}

For $T=\mathrm{PSU}_3(8)$ and $[N:\mathrm{Soc}(N)]\neq 1,2,6$, the output of \texttt{conclusion} is all \texttt{true} and this yields case (5) of Theorem \ref{new thm} except when $N \simeq \mathrm{PSU}_3(8)$, which is already covered by Theorem \ref{old thm} so this case may be ignored. Let us remark that the output of \texttt{conclusion} is \texttt{inconclusive} in this case and \cite{Tsang last} dealt with it using a slightly different method.

\vspace{1.5mm}

For $T=\mathrm{PSU}_3(8)$ and $[N:\mathrm{Soc}(N)]=2,6$, the output is as follows:
\begin{verbatim}
<2, "not normal", "false">
<6, "normal", "true">
<6, "not normal", "false">
\end{verbatim}

\noindent We then deduce case (6) of Theorem \ref{new thm}.



\section{Proof of case (1)}

Let $q = p^f$ be a prime power with $q\neq 2,3$ and let
\[ \PSL_2(q) \leq N \leq \Aut(\PSL_2(q)).\]
Since $\mathrm{Out}(\PSL_2(q))$ is abelian, such an $N$ is normal in $\Aut(\PSL_2(q))$,  so we may identify $\Aut(N)$ with $\Aut(\PSL_2(q))$, and in particular $\Inn(N)$ with $N$. By Proposition \ref{key prop}(b), it is then enough to show that there is a subgroup
\[ N \leq P \leq \Aut(\PSL_2(q))\]
for which there exist solvable subgroups $A,\, B$ of $P$ satisfying
\begin{equation}\label{want}
P = AB,\,\ A\cap B=1,\,\ AN =BN,\,\ \mbox{and $A$ splits over $A\cap N$}.
\end{equation}
Note that the last two conditions are trivial when $P$ is taken to be $N$ since then $AN = N = BN$ and $A$ obviously splits over $A\cap N = A$.

\vspace{1.5mm}

We may regard $\PSL_2(q)$ as a normal subgroup of $\PGL_2(q)$ via the natural embedding, and in particular $\PGL_2(q)$ as a subgroup of $\Aut(\PSL_2(q))$ that acts on $\PSL_2(q)$ via conjugation. It is well-known that
\begin{equation}\label{AutPSL} \Aut(\PSL_2(q)) = \PGL_2(q)\rtimes F,\end{equation}
where $F\simeq C_f$ is the subgroup consisting of the automorphisms induced by the field automorphisms of $\mathbb{F}_{p^f}$ over $\mathbb{F}_p$. By the proof of \cite[Theorem 1.3]{Tsang last}, we already know that there are solvable subgroups $C,\, D$ of $\PGL_2(q)$ such that
\[\PGL_2(q) =  CD,\,\ C\cap D= 1,\,\ C\cdot \PSL_2(q) = D\cdot \PSL_2(q).\]
Explicitly, we may take $C= \widetilde{C}/Z$ and $D = \widetilde{D}/Z$, where
\begin{align*}
 \widetilde{C} &= \left\{\begin{bmatrix} x & - dy \\ y & x-cy \end{bmatrix} : x,y\in \mathbb{F}_q,\, (x,y)\neq (0,0) \right\},\\
 \widetilde{D} & =  \left\{ \begin{bmatrix} u &v \\ 0 & w\end{bmatrix} : u,w\in\mathbb{F}_q^\times,\, v\in \mathbb{F}_q \right\}.
 \end{align*}
 Here $Z$ denotes the center of $\mathrm{GL}_2(q)$ and $X^2 + cX +d$ is the minimal polynomial  over $\mathbb{F}_q$ of some generator of $\mathbb{F}_{q^2}^\times$. Observe that
 \[ \theta \cdot\begin{bmatrix} x & y \\ z & w \end{bmatrix} Z \cdot \theta^{-1} = \begin{bmatrix}\theta(x) & \theta(y) \\ \theta(z) & \theta(w) \end{bmatrix}Z\]
 in the group (\ref{AutPSL}) for all $\theta\in F$ and $\left[\begin{smallmatrix}x & y \\ z & w\end{smallmatrix}\right]\in \mathrm{GL}_2(q)$. It in particular implies that $D$ is normalized by $F$. Thus, we have that $DE = D\rtimes E$ is a subgroup of $\Aut(\PSL_2(q))$ for any subgroup $E$ of $F$. In what follows, let us take $E$ to be the projection of $N$ onto $F$ along $\PGL_2(q)$, and let
 \[ P = \PGL_2(q)\rtimes E.\]
 Clearly $P$ contains $N$. For the choices of $A,\, B$, we separate into cases.
 
 \vspace{1.5mm}

First suppose that $q\equiv0\pmod{4}$. Then $\PGL_2(q) = \mathrm{PSL}_2(q)$ and so 
\[ N = \PGL_2(q)\rtimes E=P.\]
Take $(A,B) = (C,D\rtimes E)$. That $A,\, B$ are solvable, $P = AB$, and $A\cap B = 1$ are all clear. The last two conditions in (\ref{want}) are trivial here since $P = N$.

\vspace{1.5mm}

Next suppose that $q\equiv 1,3\pmod{4}$. It was shown in \cite[Section 4]{Tsang last} that
\begin{equation}\label{split}
 \begin{cases} C\mbox{ splits over }C\cap \PSL_2(q)&\mbox{when }q\equiv 1\hspace{-3mm}\pmod{4},\\
D\mbox{ splits over }D\cap \PSL_2(q)&\mbox{when }q\equiv 3\hspace{-3mm}\pmod{4}.
 \end{cases}\end{equation}
Here we shall take
\[ (A,B) = \begin{cases} (C,D\rtimes E)&\mbox{when }q\equiv 1\hspace{-3mm}\pmod{4},\\
(D\rtimes E,C)&\mbox{when }q\equiv 3\hspace{-3mm}\pmod{4}.
 \end{cases}\]
That $A,\, B$ are solvable, $P = AB$, and $A\cap B=1$ are all obvious. Since  
 \[ C\cdot \mathrm{PSL}_2(q) = D\cdot \mathrm{PSL}_2(q)\mbox{ and }\PSL_2(q)\leq  N,\]
we clearly have $CN = DN$. Note that $\PGL_2(q) = \PSL_2(q)D$ holds because $D$ contains an element lying outside of $\PSL_2(q)$, for example
\[ \begin{bmatrix} u & 0 \\ 0 & 1 \end{bmatrix}Z,\mbox{ where }u \in \mathbb{F}_q^\times \setminus (\mathbb{F}_q^\times)^2.\]
By definition, elements of $E$ may be represented as $\pi\eta$ for $\pi \in \PGL_2(q)$ and $\eta\in N$. Since $N$ contains $\PSL_2(q)$, we then deduce that $E$ lies in
\[ \PGL_2(q)N = \mathrm{PSL}_2(q)DN = DN.\]
Thus, we have $ (D\rtimes E)N=DN=CN$ and this yields $AN = BN$. Finally, we check that $A$ indeed splits over $A\cap N$.
\begin{itemize}
\item $q\equiv1\hspace{-1mm}\pmod{4}:$ Here $A = C$. Since $A\leq \PGL_2(q)$, we have
\[A\cap N =\begin{cases}
 A\cap \PSL_2(q) &\mbox{if }N\cap \PGL_2(q) = \PSL_2(q),\\
A & \mbox{if }N\cap \PGL_2(q) = \PGL_2(q).
\end{cases}\]
In the former case $A$ splits over $A\cap N$ by (\ref{split}), and in the latter case $A$ certainly splits over $A\cap N = A$.
\vspace{1.5mm}
\item $q\equiv3\hspace{-1mm}\pmod{4}:$ Here $A = D\rtimes E$. Since $q=p^f$, the exponent $f$ must be odd in this case, and so
\[\mathrm{Out}(\PSL_2(q))\simeq C_2\times C_f\simeq C_{2f}\]
 is cyclic. We then see that there are only two possibilities for $N$, namely
\[ N = \mathrm{PSL}_2(q)\rtimes E\mbox{ or } N = \PGL_2(q)\rtimes E.\]
In the latter case $N=P$ and $A$ certainly splits over $A\cap N = A$. In the former case, we observe that
\begin{align*}
 A\cap N& = (D\rtimes E)\cap (\PSL_2(q)\rtimes E)
=  (D\cap \PSL_2(q))\rtimes E
 \end{align*}
and in particular 
\begin{align*}
 [A:A\cap N] 
 & = [D:D\cap \PSL_2(q)].
 \end{align*}
Here, instead of applying (\ref{split}) directly, we shall need to use its proof given in \cite[Section 4]{Tsang last}. It was shown by noting that
\[  [D:D\cap \PSL_2(q)]=2,\,\  |D\cap \PSL_2(q)| =\frac{q(q-1)}{2}\]
and by using the Schur-Zassenhaus theorem. But then we have
\[ [A:A\cap N] = 2,\,\ |A\cap N| = \frac{q(q-1)}{2}\cdot |E|.\]
Since $|E|$, being a divisor of $f$, is odd in this case, we can similarly use the Schur-Zassenhaus theorem to deduce that $A$ splits over $A\cap N$.
\end{itemize}
In both cases, we have shown that $A$ splits over $A\cap N$ for our choice of $A$.

\vspace{1.5mm}

We have thus verified the existence of $P$ satisfying (\ref{want}) in all cases. This proves case (2) of Theorem \ref{new thm}.


\end{document}